\documentclass[11pt, reqno]{amsart}

\usepackage[OT2, T1]{fontenc}
\usepackage{url}
\usepackage{amsmath}
\usepackage{array}
\usepackage{graphicx}
\usepackage{amsfonts}
\usepackage{amssymb}
\usepackage{amsthm}
\usepackage{hyperref}
\usepackage{colonequals}	% for nice :=
\usepackage{color}
\usepackage{mathabx}		% for \widec
\usepackage{mathrsfs}
\usepackage{amscd}
\usepackage[all,cmtip]{xy}	% for commutative diagrams
\usepackage{sseq}			% for spectral sequences
\usepackage{verbatim}
\usepackage{parskip}
\usepackage[in]{fullpage}
\usepackage{mathrsfs} 			% for \mathscr (script letters)
\usepackage{bm} 				% for bold Greek letters
\usepackage{cleveref}
\usepackage{enumitem}
\usepackage{moreenum}			% for enumeration via Greek letters
\usepackage[alphabetic,lite]{amsrefs} 	% for bibliography
\usepackage{stmaryrd}	% this is for \llbracket and \rrbracket

\DeclareSymbolFont{cyrletters}{OT2}{wncyr}{m}{n}
\DeclareMathSymbol{\Sha}{\mathalpha}{cyrletters}{"58}

\usepackage{color}

\newcommand{\bC}{\mathbb{C}}

\newcommand{\bF}{\mathbb{F}}

\newcommand{\bQ}{\mathbb{Q}}

\newcommand{\bZ}{\mathbb{Z}}

\newcommand{\cO}{\mathcal{O}}

\newcommand{\fm}{\mathfrak{m}}

\newcommand{\fp}{\mathfrak{p}}

\newcommand{\ra}{\rightarrow}

\newcommand{\xra}{\xrightarrow}

\newcommand{\hra}{\hookrightarrow}

\newcommand{\I}{^{\infty}}
\newcommand{\wt}{\widetilde}
\newcommand{\wh}{\widehat}
\newcommand{\eps}{\epsilon}

\newcommand{\pr}{^{\prime}}
\newcommand{\prpr}{^{\prime\prime}}
\newcommand{\ce}{\colonequals}
\newcommand{\ov}{\overline}

\renewcommand{\b}{\textbf}
\newcommand{\surjects}{\twoheadrightarrow}

\newcommand{\tensor}{\otimes} 		% binary tensor product
\newcommand{\isomto}{\overset{\sim}{\longrightarrow}}

\newcommand{\nr}{{\mathrm{nr}}}		% unramified cohomology or max unramified extension
\newcommand{\normal}{{\lhd}}			% normal subgroup
\newcommand{\llb}{\llbracket}		% [[
\newcommand{\rrb}{\rrbracket}		% ]]
\renewcommand{\i}{^{-1}}
\renewcommand{\th}{^{\mathrm{th}}}

\providecommand{\abs}[1]{\left\lvert#1\right\rvert}

\providecommand{\p}[1]{\left(#1\right)}

\providecommand{\f}[2]{\frac{#1}{#2}}
\DeclareMathOperator{\Ker}{Ker}			% Kernel
\DeclareMathOperator{\Spec}{Spec}		% Spectrum of a ring
\DeclareMathOperator{\Hom}{Hom}			% Set of arrows between two object
\DeclareMathOperator{\Char}{char}		% Characteristic of a field
\DeclareMathOperator{\Frac}{Frac}		% Field of fractions
\DeclareMathOperator{\Ht}{ht}			% Height of an ideal
\DeclareMathOperator{\Gal}{Gal}	% Galois group
\DeclareMathOperator{\Tr}{Tr}		% Trace
\DeclareMathOperator{\ab}{ab}		% abelianization
\DeclareMathOperator{\Ind}{Ind}		% Induced representation
\DeclareMathOperator{\Res}{Res}		% Restriction of the representation
\DeclareMathOperator{\GL}{GL}		% The general linear group
\DeclareMathOperator{\End}{End}		% The algebra of endomorphisms
\DeclareMathOperator{\Aut}{Aut}		% The group of automorphisms
\DeclareMathOperator{\rk}{rk}		% rank
\DeclareMathOperator{\Art}{Art}		% local Artin homomorphism
\DeclareMathOperator{\Sw}{Sw}		% Swan conductor
\DeclareMathOperator{\Frob}{Frob}		% Frobenius
\newcommand{\ba}{\begin{aligned}}
\newcommand{\ea}{\end{aligned}}
\newcommand{\be}{\begin{equation}}
\newcommand{\ee}{\end{equation}}
\newcommand{\pf}{\begin{proof}}
\newcommand{\bpf}{\begin{proof}}
\newcommand{\epf}{\end{proof}}
\newcommand{\bthm}{\begin{thm}}
\newcommand{\ethm}{\end{thm}}
\newcommand{\bprop}{\begin{prop}}
\newcommand{\eprop}{\end{prop}}
\newcommand{\bcor}{\begin{cor}}
\newcommand{\ecor}{\end{cor}}
\newcommand{\brem}{\begin{rem}}
\newcommand{\erem}{\end{rem}}
\newcommand{\brems}{\begin{rems} \hfill \begin{enumerate}[label=\b{\thesubsection.},ref=\thesubsection]}
\newcommand{\remi}{\addtocounter{subsection}{1} \item}
\newcommand{\erems}{\end{enumerate} \end{rems}}
\newcommand{\blem}{\begin{lemma}}
\newcommand{\elem}{\end{lemma}}
\newcommand{\bconj}{\begin{conj}}
\newcommand{\econj}{\end{conj}}
\newcommand{\bprob}{\begin{Problem}}
\newcommand{\eprob}{\end{Problem}}
\newcommand{\bq}{\begin{q}}
\newcommand{\eq}{\end{q}}
\newcommand{\benum}{\begin{enumerate}[label={(\alph*)}]}
\newcommand{\eenum}{\end{enumerate}}
\newcommand{\bc}{\begin{comment}}
\newcommand{\ec}{\end{comment}}
\newcommand{\bcl}{\begin{claim}}
\newcommand{\ecl}{\end{claim}}
\newcommand{\beg}{\begin{eg}}
\newcommand{\eeg}{\end{eg}}
\newcommand{\lab}{\label}

\theoremstyle{plain}
\newtheorem{thm}[subsection]{Theorem}
\Crefname{thm}{Theorem}{Theorems}

\Crefname{rethm}{Theorem}{Theorem}
\newtheorem{prop}[subsection]{Proposition}
\Crefname{prop}{Proposition}{Propositions}
\newtheorem{q}[subsection]{Question}
\Crefname{q}{Question}{Questions}
\Crefname{eg}{Example}{Examples}
\newtheorem{Problem}[subsection]{Problem}
\Crefname{Problem}{Problem}{Problems}
\newtheorem{conj}[subsection]{Conjecture}
\Crefname{conj}{Conjecture}{Conjectures}
\newtheorem{cor}[subsection]{Corollary}
\Crefname{cor}{Corollary}{Corollaries}

\newtheorem{lemma}[subsection]{Lemma}

\newtheorem{lem}[equation]{Lemma}
\Crefname{lem}{Lemma}{Lemmas}

\theoremstyle{remark}
\newtheorem{claim}[equation]{Claim}
\Crefname{claim}{Claim}{Claims}

\theoremstyle{definition}
\newtheorem{eg}[subsection]{Example}
\newtheorem{rem}[subsection]{Remark}
\Crefname{rem}{Remark}{Remarks}
\newtheorem*{rems}{Remarks}

\newtheoremstyle{subsection-tweak}
   {11pt}
   {3pt}%
   {}
   {}%
   {\bfseries}
   {}%
   {.5em}
   {\thmnumber{\@{#1}{}\@{#2}.}%
    \thmnote{~{\bfseries#3.}}}

\Crefname{innercustomconj}{Conjecture}{Conjecture}

\theoremstyle{subsection-tweak}
\newtheorem{pp}[subsection]{}
\newcommand{\bpp}{\begin{pp}}
\newcommand{\epp}{\end{pp}}

\numberwithin{equation}{subsection}

\begin{document}
\author{K\k{e}stutis \v{C}esnavi\v{c}ius}
\title{Local factors valued in normal domains}
\date{\today}
\address{Department of Mathematics, University of California, Berkeley, CA 94720-3840, USA}
\email{kestutis@berkeley.edu}
\urladdr{http://math.berkeley.edu/~kestutis/}

\begin{abstract}
We give an exposition of Deligne's theory of local $\eps_0$-factors over fields and discrete valuation rings under the assumption that the theory over the complex numbers is known. We then employ standard techniques from algebraic geometry to deduce the theory of local $\eps_0$-factors over arbitrary normal integral schemes.
\end{abstract}

\maketitle

\section{Introduction} \lab{intro}

In \cite{Del73}*{\S4}, Deligne presented an elegant argument proving the existence of the theory of $\eps$-factors of local Weil representations over the complex numbers; this theory had previously been predicted by Langlands. For Weil representations of a nonarchimedean local field $K$ of residue characteristic $p$, Deligne proceeded to show in \cite{Del73}*{\S6} how to use the established $F = \bC$ case to deduce a similar ``mod $l$'' theory over every field $F$ of characteristic $l$ different from $p$. Both the complex and the mod $l$ theories have subsequently been of significant importance, for instance, in considerations concerning the local Langlands correspondence.

Deligne's argument concerning the mod $l$ theory is brief, with many details left to the reader. Our goal here is to recall it in full and combine its ideas with standard techniques from algebraic geometry to deduce the theory of local $\eps_0$-factors not only over fields but also over arbitrary normal integral schemes on which $p$ is invertible. In precise terms, assuming the theory of \cite{Del73}*{\S4} over the complex numbers as known, we prove

\bthm \lab{main}
There is a unique assignment $\epsilon_0$ which to the data of 
\begin{itemize}
\item 
A nonarchimedean local field $K$ with the ring of integers $\cO_K$ and a finite residue field of characteristic $p$,

\item A separable closure $K^s$ of $K$,

\item A normal integral $\bZ[\f{1}{p}]$-scheme $S$,

\item A continuous representation $V$ of the Weil group $W(K^s/K)$ over $S$ (cf.~\S\ref{weil} and \S\ref{Gr-rings}),

\item A nontrivial additive character $\psi\colon (K, +) \ra \Gamma(S, \cO_S^\times)$ (cf.~\S\ref{add-char}), and

\item A $\Gamma(S, \cO_{S})$-valued Haar measure $C \mapsto \int_C dx$ on $K$ such that $\int_{\cO_K} dx \in \Gamma(S, \cO_S^\times)$ (cf.~\S\ref{haar})
\end{itemize}

associates $\epsilon_0(V, \psi, dx) \in \Gamma(S, \cO_S^\times)$ in such a way that

\begin{enumerate}[label={(\roman*)}] 
\item \lab{BC}
The formation of $\epsilon_0(V, \psi, dx)$ is compatible with base change: for a morphism $f: S\pr \ra S$ of normal integral $\bZ[\f{1}{p}]$-schemes and an $S$-representation $V$, abusing the $f^*$ notation one has
\[
\eps_0(f^*V, f^{*}\psi, f^{*}dx) = f^{*}(\eps_0(V, \psi, dx)).
\]

\item \lab{add}
Exactness of a sequence $0 \ra V\pr \ra V \ra V\prpr \ra 0$ entails
\[
\eps_0(V, \psi, dx) = \eps_0(V\pr, \psi, dx)\eps_0(V\prpr, \psi, dx).
\]
In particular, $\epsilon_0(V, \psi, dx)$ depends only on the class of $V$ in the Grothendieck group $R_S(W(K^s/K))$ and makes sense for every $v \in R_S(W(K^s/K))$ (see \S\ref{Gr-rings} for the definition of $R_S(W(K^s/K))$).

\item \lab{haar-scale}
For $a\in \Gamma(S, \cO_S^\times)$, one has 
\[
\eps_0(V, \psi, a\cdot dx) = a^{\rk V} \eps_0(V, \psi, dx).
\]

\item \lab{ind}
For a finite subextension $K^s/L/K$ and a virtual representation $v$ of $W(K^s/L)$ of rank $0$,
\[ 
\eps_0\p{\Ind_{W(K^s/L)}^{W(K^s/K)} v, \psi} = \eps_0(v, \psi \circ \Tr_{L/K}),
\]
where $\Ind_{W(K^s/L)}^{W(K^s/K)} v \ce \cO_S[W(K^s/K)] \tensor_{\cO_S[W(K^s/L)]} v$. (Omitting $dx$ is justified by \ref{haar-scale}.)

\item \lab{dim-1}
If $V$ is of dimension $1$ (i.e., a line bundle) and $\chi\colon W(K^s/K) \ra \Gamma(S, \cO_S^\times)$ is the character giving the Weil group action, then 
\be\lab{dim-1-def}\tag{$\bigstar$}
\eps_0(\chi, \psi, dx) = \int_{\gamma\i \cO_K^\times} \chi\i(x)\psi(x) dx
\ee
where $\gamma \in K^\times$ is an element of valuation $\Sw(\chi) + n(\psi) + 1$ (for the definition of the Swan conductor $\Sw(\chi) \in \bZ_{\ge 0}$, that of $n(\psi) \in \bZ$, and the meaning of the integral, see \Cref{Swan-sch}, \S\ref{add-char}, and \S\ref{haar}).
\eenum
\ethm

\brems
\remi \lab{hist}
Restricting to the main case of interest, $S = \Spec \bC$, the existence and uniqueness of an $\eps_0$ satisfying \ref{add}--\ref{dim-1} was envisioned by Langlands: the automorphic side of his conjectural correspondence features decompositions of signs of functional equations of global $L$-functions as products of local factors, and \ref{add}--\ref{dim-1} are the defining properties of the corresponding local factor on the Galois side. In \cite{Lan70}, Langlands attempted to give a local proof of the existence of $\eps_0$ but abandoned the project once Deligne found a short proof \cite{Del73}*{\S4}, which, however, uses global arguments. Publishing a complete local proof remains an outstanding problem. For further and more accurate historical remarks, see the website of \cite{Lan70}.

\remi \lab{S=C}
For the purpose of proving \Cref{main}, we will take the existence and uniqueness of an $\eps_0$ satisfying \ref{add}--\ref{dim-1} with $S = \Spec \bC$ as known. Although existence is intricate, uniqueness follows readily from a suitable version of Brauer's induction theorem due to the imposed \eqref{dim-1-def} in the $1$-dimensional case: see \cite{Tat79}*{2.3.1} or the proof of \Cref{alg-cl} \ref{unique} below. Tate's thesis \cite{Tat50}, which has been a major influence for the ideas mentioned in Remark \ref{hist}, motivates \eqref{dim-1-def} and also provides the main input for proving that \eqref{dim-1-def} results in a global unit on $S$, as is implicit in \Cref{main} and will be argued in the course of the proof.

\remi \lab{full-eps}
When $S = \Spec F$ for a field $F$, e.g., when $S =\Spec \bC$, it is commonplace to work with 
\[
\eps(V, \psi, dx) \ce \eps_0(V, \psi, dx)\det(-\mathrm{Frob}_K\, |\, V^{I_K})\i,
\]
which is the \emph{local $\eps$-factor} of $V$ (for the choice of $\psi$ and $dx$); here $\Frob_K$ is a geometric Frobenius defined in \S\ref{not}. However, \Cref{main} concerns $\eps_0$ instead because the formation of $\eps$ is not compatible with base change, i.e., the analogue of \ref{BC} fails for $\eps$.

\remi
We collect formulas concerning $\eps_0$ in \S\ref{formulas}. The only one of these that exhibits phenomena not observed when $S = \Spec \bC$ is that for the inverse in the higher-dimensional case, see \S\ref{high-dim}.

\remi 
The possibility of extending Deligne's theory of $\eps_0$-factors for representations of $W(K^s/K)$ over fields to those over more general coefficient rings $R$ has also been considered in \cite{Yas09}, where such an extension is proposed for Noetherian local rings $R$ that have an algebraically closed residue field of characteristic different from $p$ and satisfy $R^{\times p} = R^\times$. In the case of an algebraically closed field of characteristic different from $p$, the $\eps_0$-factors of op.~cit.~agree with those of \Cref{main} due to \cite{Yas09}*{Thm.~1.1 (3)}. For general $R$ as above that are also normal domains, if the $\eps_0$-factors there are compatible with base change along not necessarily local homomorphisms $R \ra R\pr$, then one gets the similar agreement by taking an algebraic closure of $\Frac R$ for $R\pr$.
\erems

%\kestutis{Doing reduced schemes rather than integral schemes seems difficult: although for a reduced ring $A$ one has the injection $A \hra \prod_{\fp} A/\fp$ where the product runs over the min'l primes of $A$, it is not clear that the resulting product of $\eps_0$'s lies in the image of $A$. Especially troubling is the example case $A = k[x, y]/(y(x^2 - y))$ with $\fp_1 = (y)$, $\fp_2 = (x^2 - y)$ when using the exact sequence $0 \ra A/\fp_1\cap \fp_2 \ra A/\fp_1 \oplus A/\fp_2 \ra A/(\fp_1 + \fp_2) \ra 0$ of $A$-modules, to check that $\eps_0$ is in $A/\fp_1 \cap \fp_2 = A$, one would need to check that the images of the two $\eps_0$'s in $A/(\fp_1+ \fp_2)$ agree. However, the latter is $k[x]/(x^2)$, which is nonreduced, and hence I cannot use my understanding of $\eps_0$ in the integral case to say anything (a priori) about $\eps_0$ over $k[x]/(x^2)$.}

\bpp[Notation] \lab{not}
For a field $F$, its fixed choices of separable and algebraic closures are denoted by $F^s \subset \ov{F}$. As mentioned above, $K$ is a nonarchimedean local field, whereas $\cO_K$ and $\bF_K$ are its ring of integers and residue field (so $\bF_K$ is a finite field); $p = \Char \bF_K$ is the residue characteristic. A \emph{geometric Frobenius} is any $\Frob_K \in W(K^s/K) \subset \Gal(K^s/K)$ whose inverse reduces to the Frobenius automorphism $x \mapsto x^{\#\bF_K}$ in $\Gal(\ov{\bF}_K/\bF_K)$.  For an integer $n \ge 1$, a primitive $n\th$ root of unity is denoted by $\zeta_n$. For a scheme $S$, the local ring and the residue field of a point $s \in S$ are denoted by $\cO_{S, s}$ and $k(s)$.
 \epp

\bpp[Conventions]
The following assumptions are implicit throughout: all rings are commutative and unital; all representations are of finite rank; all representations of $I_K$ or $W(K^s/K)$ are trivial on an open subgroup of $I_K$.
\epp

\subsection*{Acknowledgements}
I thank Pierre Deligne for explaining me the argument used to prove \Cref{fields}. I thank Bjorn Poonen, Jack Thorne, Seidai Yasuda, and the referee for helpful comments and suggestions. This note was written as a term paper for the course Topics in Automorphic Forms that the author took in Fall 2013 at Harvard University. I thank Jack Thorne for an interesting and useful course and for an opportunity to write up the results presented here.

\section{Recollections} \lab{defn}

We gather several relevant concepts and constructions that are used freely in other sections.

\bpp[The Weil group $W(K^s/K)$ (see also \cite{BH06}*{\S28} if needed)] \lab{weil}
It is the subgroup of those elements of $\Gal(K^s/K)$ that reduce to an integral power of Frobenius in $\Gal(\ov{\bF}_K/\bF_K)$. Thus,
\[
1 \ra I_K \ra W(K^s/K) \ra \bZ \ra 1
\]
is short exact, where $I_K \normal \Gal(K^s/K)$ is the inertia. The Weil group is topologized by insisting that $I_K$ with its profinite topology be an open subgroup. 

The open subgroups of $W(K^s/K)$ of finite index are precisely the $W(K^s/L)$ for finite subextensions $K^s/L/K$; normality of the subgroup corresponds to $L/K$ being Galois. The natural map 
\[
W(K^s/K)/W(K^s/L) \ra \Gal(K^s/K)/\Gal(K^s/L)
\]
is bijective. In particular, $W(K^s/K) \subset \Gal(K^s/K)$ is dense and the representations of $W(K^s/K)$ whose kernel is open and of finite index are identified with those of $\Gal(K^s/K)$. Such representations of $W(K^s/K)$ are called \emph{Galois}.

Let $W(K^s/K)^{\ab}$ be the maximal abelian Hausdorff quotient of $W(K^s/K)$. Due to the ambiguity in choosing $K^s$, the Weil group is determined by $K$ only up to an inner automorphism of $\Gal(K^s/K)$; this ambiguity disappears for $W(K^s/K)^{\ab}$, which is determined by $K$ up to a unique isomorphism. Local class field theory furnishes the local Artin homomorphism 
\be\lab{Art}
\Art_K\colon K^\times \ra W(K^s/K)^{\ab},
\ee
which is an isomorphism of topological groups. We choose to normalize it so that uniformizers are brought to geometric Frobenii. As usual, we identify continuous $1$-dimensional characters of $W(K^s/K)$ valued in Hausdorff abelian groups with those of $K^\times$ by means of $\Art_K$. We let $\abs{\cdot}_K\colon K^\times \ra \bZ[\f{1}{p}]^\times$ be the unramified character that takes the value $(\#\bF_K)\i$ on uniformizers.
\epp

\bpp[Grothendieck groups] \lab{Gr-rings}
A \emph{representation} of a discrete group $G$ over a scheme $S$ is an $\cO_S$-module $V$ that is locally free of finite rank and is endowed with an $\cO_S$-linear action of $G$. For the topological groups $I_K$ and $W(K^s/K)$, one only considers \emph{continuous} representations, i.e., those $V$ on which an open subgroup of $I_K$ acts trivially (some authors call such representations \emph{smooth}).

Let $R_S(G)$ (resp.,~$R_S(I_K)$ or $R_S(W(K^s/K))$) be the \emph{Grothendieck group} of representations of $G$ (resp.,~continuous representations of $I_K$ or $W(K^s/K)$) over $S$, i.e.,  $R_S(G)$ is the quotient of the free abelian group on the set $\{[V]\}$ of isomorphism classes of representations of $G$ over $S$ by the subgroup generated by the relations $[V] = [V\pr] + [V^{\prime\prime}]$ for all exact sequences
\[
0 \ra V\pr \ra V \ra V^{\prime\prime} \ra 0,
\]
and likewise for $I_K$ or $W(K^s/K)$ in place of $G$. A \emph{virtual representation} is an element of a Grothendieck group. The rank of $V$ is an integer valued locally constant function on $S$; its additivity in short exact sequences defines the notion of the \emph{rank} of a virtual representation.

Letting $J$ run over the open subgroups of $I_K$ that are normal in\footnote{Every open $J\pr\normal I_K$ contains such a $J$ of the form $\bigcap_{i = 0}^{n - 1} (\Frob_K^i J\pr \Frob_K^{-i})$ for some $n \ge 1$: indeed, $J\pr$ corresponds to a finite Galois $L/K^\nr$, which descends to a Galois extension of the degree $n$ unramified extension of $K$ for some~$n \ge 1$.} $W(K^s/K)$, one has
\be\lab{R-I-W}
R_S(I_K) = \textstyle{\varinjlim_J} R_S(I_K/J) \quad\quad\text{and}\quad\quad R_S(W(K^s/K)) = \varinjlim_{J} R_S(W(K^s/K)/J).
\ee
Due to the $\cO_S$-flatness of $V$, the tensor product endows $R_S(-)$ with the structure of a commutative ring with $[\cO_S]$ as the multiplicative unit. The subset $R_S^0(-) \subset R_S(-)$ of virtual representations that are of rank $0$ at every $s\in S$ is an ideal. Pullback along $S\pr \ra S$ induces a ring homomorphism $R_S(-) \ra R_{S\pr}(-)$, which maps $R_S^0(-)$ to $R_{S\pr}^0(-)$. If $S = \Spec A$, one often writes $R_A(-)$ for $R_S(-)$.
\epp

\bpp[Additive characters]\lab{add-char}
For an integral $\bZ[\f{1}{p}]$-scheme $S$, an \emph{additive character}
\[
\psi\colon (K, +) \ra \Gamma(S, \cO_S^\times)
\]
is a locally constant abelian group homomorphism as indicated. Local constancy is equivalent to the existence of an integer $n$ such that $\psi|_{\pi^{-n}\cO_K} = 1$, where $\pi \in \cO_K$ is a uniformizer, and forces the values of $\psi$ to be $p$-power roots of unity (if $\Char K = p$, these values are even $p\th$ roots of unity). If $\psi$ is nontrivial, as we assume from now on, we let $n(\psi)$ be the maximal $n$ as above. Since $p$ is a unit on $S$, a $p$-power root of unity in $\Gamma(S, \cO_S^\times)$ has trivial image in $k(s)^\times$ for an $s\in S$ if and only if it is trivial to begin with. Consequently, $n(\psi) = n(\psi_{k(s)})$ for every $s \in S$, and hence $n(\psi)$ is stable under pullback along every $f\colon S\pr \ra S$ with an integral $S\pr$.

The group $K^\times$ acts freely on the set of nontrivial $\psi$ by setting $(a\psi)(x) \ce \psi(ax)$ for $a \in K^\times$. If the set is nonempty, then the action is also transitive: endowing $S$ with a structure of an $R$-scheme, where $R = \bZ[\zeta_{p\I}]$ if $\Char K = 0$ and $R = \bZ[\zeta_p]$ if $\Char K = p$, realizes every $\psi$ as the pullback of an additive character valued in $R^\times$, and hence reduces the transitivity claim to the classical $S = \Spec \bC$ case treated, e.g., in \cite{BH06}*{\S1.7, Prop.}. In conclusion, if the set of nontrivial $\psi$ is nonempty, then it has a natural structure of a $K^\times$-torsor.
\epp

\bpp[Haar measures valued in abelian groups (\cite{Del73}*{6.1})] \lab{haar}
Fix a $\bZ[\f{1}{p}]$-module $A$ and let $C$ range over the compact open subsets of $K$. An $A$-valued \emph{Haar measure} on $K$ is a function $C \mapsto \int_C dx \in A$ that is translation invariant and additive in disjoint unions. Since every $C$ is a disjoint union of translates of balls centered at $0 \in K$ and $A$ is uniquely $p$-divisible, a choice of an $A$-valued Haar measure amounts to that of the element $\int_{\cO_K} dx \in A$. Once the choice is made, if $A$ is in addition a ring, one can integrate locally constant compactly supported $f\colon K \ra A$ and write $\int_K f(x) dx$ for the resulting finite sums; if $f$ is implicitly multiplied by the characteristic function of $C$, one writes $\int_C f(x) dx$ instead. Coupled with \eqref{Art}, this clarifies \eqref{dim-1-def}, which takes $A = \Gamma(S, \cO_S)$.
\epp

\bpp[Sufficiently large fields] \lab{suff-large}
Fix a finite group $G$ and let $m$ be the least common multiple of the orders of elements of $G$. A field $F$ is \emph{sufficiently large} if it contains the $m\th$ roots of unity. A separably closed $F$ is sufficiently large for every $G$. 
\epp

\bpp[Virtual representations over fields] \lab{Gr-fields}
For a field $F$, thanks to the Jordan--H\"{o}lder theorem~for abelian categories \cite{Ses67}*{Thm.~2.1}, $R_F(G)$ of \S\ref{Gr-rings} is the free abelian group on the set of isomorphism classes of irreducible representations of $G$ over $F$, and similarly for $I_K$ or $W(K^s/K)$ and representations that are trivial on an open subgroup of $I_K$ (cf.~\cite{Ser77}*{\S14.1, Prop.~40} if needed). 

If $V$ and $\wt{V}$ are nonisomorphic irreducible representations of a finite group $G$, then the extensions of scalars $V_{F\pr}$ and $\wt{V}_{F\pr}$ to every overfield $F\pr/F$ have no common composition factors \cite{CR81}*{Ex.~7.9}. Consequently, $R_F(G) \ra R_{F\pr}(G)$ is injective. It is also surjective if $F$ is sufficiently large, as is clear from Brauer's induction theorem \ref{brau} \ref{brau-a}, whose proof does not use this surjectivity. Therefore, for sufficiently large $F$, extension of scalars induces a bijection between the sets of isomorphism classes of irreducible representations of $G$ over $F$ and those over $F\pr$. If $F$ contains the $m\th$ roots of unity for every $m$, e.g., if $F$ is separably closed, then this bijection remains in place for $I_K$, because $R_F(I_K) \ra R_{F\pr}(I_K)$ is an isomorphism thanks to the previous discussion and \eqref{R-I-W}.
\epp

\bpp[The decomposition homomorphism] \lab{dec}
Fix a finite group $G$ and a discrete valuation ring $A$ with the fraction field $\eta$ of characteristic $0$ and the residue field $F = A/\fm$ of characteristic $l$. For a representation $V$ of $G$ over $\eta$, a choice of a $G$-stable $A$-lattice $\Lambda \subset V$ gives rise to the representation $\Lambda/\fm \Lambda$ of $G$ over $F$ whose class in $R_F(G)$ does not depend on $\Lambda$ \cite{Ser77}*{\S15.2, Thm.~32 and Rem.~(1)} (the $l > 0$ assumption of loc.~cit.~is not used for this). The resulting \emph{decomposition homomorphism}
\[
d_G\colon R_\eta(G) \ra R_F(G)
\]
preserves ranks and commutes with restriction and induction. If $A$ is complete, then $d_G$ is surjective: Cohen's structure theorem \cite{Mat89}*{28.3 (ii)} settles the equicharacteristic case, whereas \cite{Ser77}*{\S16.1, Thm.~33} treats the case $l > 0$. As noted in \cite{Ser77}*{Ex.~16.1}, surjectivity also holds if $\eta$ is sufficiently large\footnote{However, $d_G$ is not surjective in general, see \cite{Ser77}*{Ex.~16.2} or \cite{CR81}*{p.~512, Ex.~21.4}.} because $d_G$ commutes with $R_\eta(G) \ra R_{\wh{\eta}}(G)$, which is an isomorphism \cite{Ser77}*{\S12.3}. 

In the equicharacteristic $0$ case, $d_G$ is injective without additional assumptions on $A$: for every virtual character $\chi$, one has $\langle \chi, \chi \rangle \in \bZ_{\ge 0}$ ; moreover, $\langle \chi, \chi \rangle = 0$ if and only if $\chi = 0$. Consequently, if in this case $A$ is complete or $\eta$ is sufficiently large, then $d_G$ is an isomorphism.
\epp

\bpp[The Swan representation] \lab{Sw-rep}
Let $J$ be a continuous finite quotient of $I_K$. Continuity means the openness of $\Ker(I_K \surjects J)$ and is a nonvacuous condition, as we now explain in a digression. By Krasner's lemma, there are only countably many finite degree subextensions of $K^s/K$ because the same holds for global $K$. Applying this observation to finite unramified extensions of $K$, we conclude that the same holds for $K^s/K^{\nr}$, in other words, that $I_K$ has only countably many open subgroups of finite index. On the other hand, local class field theory applied to finite unramified extensions of $K$ produces continuous surjections $I_K \surjects (\bZ/p\bZ)^n$ for arbitrarily large $n \in \bZ_{> 0}$. Thus, $\bigoplus_{i = 0}^\infty \bZ/p\bZ$ is a quotient of $I_K$, so that $I_K$ has uncountably many distinct surjections onto $\bZ/p\bZ$. For cardinality reasons, one of these surjections must have a kernel that is not open. 

Let $J = J_0 \rhd J_1 \rhd \dotsb$ be the ramification filtration in the lower numbering, so the $J_i$ are normal in $J$ and $J_1$ is the image of the wild inertia. The \emph{Artin character} of $J$ is the class function
\be \lab{artin-char}
a_J \ce \sum_{i = 0}\I \f{1}{[J_0 : J_i]} \Ind_{J_i}^{J_0} u_i,
\ee
where $u_i$ is the augmentation representation of $J_i$, i.e., 
\[
u_{i}(j) = \begin{cases} - 1,\ \,\quad\quad  \text{ if $j \in J_i \setminus \{1\}$,} \\ \# J_i -1,\ \text{ if $j = 1$}, \end{cases} \quad \text{so $r_{J_i} = \b{1}_{J_i} \oplus u_i$, where $r_{J_i}$ is the regular representation.}
\]
The formula for the induced character shows that the sum in \eqref{artin-char} is finite and $a_J$ is $\bZ$-valued; moreover, as the name suggests, $a_J$ is the character of a complex representation \cite{Ser79}*{VI.\S2, Thm.~1 and Prop.~2}, namely, the \emph{Artin representation} of $J$. The \emph{Swan character}
\[
\Sw_J \ce \sum_{i = 1}\I \f{1}{[J_0 : J_i]} \Ind_{J_i}^{J_0} u_i
\]
inherits these properties: a priori $\Sw_J$ is a virtual character with $\langle \Sw_J, \chi \rangle \ge 0$ for every character $\chi$ of $J$, and hence a posteriori the character of a complex representation, namely, the \emph{Swan representation} of $J$. Moreover, $\Sw_J$ vanishes on $J \setminus J_1$, so whenever the Swan character can be realized as a representation over the fraction field $\eta$ of characteristic $0$ of a Dedekind domain $A$ with $p \in A^\times$, due to \cite{Swa63}*{Thm.~5}, it can also be realized as a finite projective $A[J]$-module. Realizability over a sufficiently large $\eta$ is automatic (cf.~\S\ref{Gr-fields}); therefore, \cite{Ser77}*{\S16.3, Prop.~44} realizes $\Sw_J$ uniquely as a finite projective $\bZ_l[J]$-module for every $l \neq p$. We continue to write $\Sw_J$ for its base change to a $\bZ_l$-algebra $A$, for instance, to a field of characteristic $l$; the resulting $\Sw_J$ is a finite projective $A[J]$-module.

If $J\pr$ is a quotient of $J$, then $a_{J\pr} \cong a_J \tensor_{\bC[J]} \bC[J\pr]$ \cite{Ser79}*{VI.\S2 Prop.~3}; the same relation holds for the augmentation representations, so also $\Sw_{J\pr} \cong \Sw_J \tensor_{\bC[J]} \bC[J\pr]$. Moreover, uniqueness of the realization of $\Sw_{J\pr}$ as a projective $\bZ_l[J\pr]$-module for $l \neq p$ entails the $A = \bZ_l$ case of the isomorphism
\be\lab{Sw-ind}
\Sw_{J\pr} \cong \Sw_J \tensor_{A[J]} A[J\pr] \quad\quad \text{for every $\bZ_l$-algebra $A$},
\ee
and the general case follows by base change to $A$.
\epp

\bpp[The Swan conductor]
Let $V$ be a continuous representation of $I_K$ over a field $F$ of characteristic $l$ with $l \neq p$, and let $J$ be a continuous finite quotient through which the $I_K$-action on $V$ factors. The \emph{Swan conductor}\footnote{Some authors call $\Sw V$ \emph{the exponent of the Swan conductor}, reserving the term \emph{Swan conductor} for the corresponding power of the maximal ideal of $\cO_K$.} of $V$ is 
\be\lab{Sw-def}
\Sw V \ce \dim_F \Hom_J(\Sw_J, V),
\ee
where in the $l > 0$ case one uses the projective $F[J]$-module $\Sw_J$ defined in \S\ref{Sw-rep}, and in the $l = 0$ case one writes $\chi$ for the character of $V$ and interprets the right hand side of \eqref{Sw-def} as
\be\lab{Sw-char-0}
\langle \Sw_J, \chi \rangle = \f{1}{\#J} \sum_{j \in J}  \Sw_J(j)\chi(j).
\ee
Of course, if $l = 0$ and $\Sw_J$ is realizable over $F$, then \eqref{Sw-def} and \eqref{Sw-char-0} agree thanks to character theory. Moreover, one can assume realizability for the purpose of the definition, because $\Sw V$ is invariant under base change to an overfield $F\pr/F$ regardless of $l$. Likewise, $\Sw V$ is invariant under change of $J$ due to \eqref{Sw-ind} and the adjunction $-\tensor_{F[J]} F[J\pr] \dashv \Hom_{J\pr}(F[J\pr], -)$. Since it is also additive in exact sequences due to projectivity of $\Sw_J$, it extends to a homomorphism~$\Sw\colon R_F(I_K) \ra \bZ$. 

To define $\Sw V$, it is not necessary to restrict to representations over fields, as \Cref{Swan-sch} below shows. For its proof, we recall a well-known \Cref{lem-EGA}, which will also be used later.
\epp

\blem\lab{lem-EGA}
For a locally Noetherian scheme $S$, a point $s \in S$, and its specialization $s\pr \neq s$, there is a complete discrete valuation ring $A$ and a morphism $\Spec A \ra S$ mapping the generic and the closed points of $\Spec A$ to $s$ and $s\pr$, respectively.
\elem

\bpf
Replace a discrete valuation ring provided by \cite{EGAII}*{7.1.9} by its completion.
\epf

\bprop\lab{Swan-sch}
Let $V$ be a continuous representation of $I_K$ over an integral $\bZ[\f{1}{p}]$-scheme $S$ (cf.~\S\ref{Gr-rings}). For varying $s\in S$, the Swan conductor of the residual representation $V_{k(s)}$ is constant. 
\eprop

The common value of the $\Sw V_{k(s)}$ is the \emph{Swan conductor} of $V$. It has already been used in \ref{dim-1}.

\bpf
Since every two nonempty opens of $S$ intersect, it suffices to treat the affine case $S = \Spec A$ and assume that $V$ is free. Writing $A$ as a filtered direct limit of finite type $\bZ[\f{1}{p}]$-subalgebras, one uses limit arguments and invariance of $\Sw V_{k(s)}$ under base change to overfields to assume further that $A$ is Noetherian. Taking $s$ in \Cref{lem-EGA} to be the generic point, one finally reduces to the case when $A$ is a complete discrete valuation ring with the fraction field $\eta$ and the residue field $F$. 

In this case, if $\Char F = 0$, then the claim follows from \eqref{Sw-char-0}, which takes values in $A$. If, on the other hand, $\Char F = l$ with $l > 0$, then $A$ is a $\bZ_l$-algebra and the projectivity of the finite $A[J]$-module $\Sw_J$ realizes it as a direct summand of a finite free $A[J]$-module; consequently, the $A$-module $\Hom_J(\Sw_J, V)$ is also finite free, and it remains to note that the rank of its pullback to $\eta$ (resp.,~$F$) equals $\Sw V_\eta$ (resp.,~$\Sw V_F$).
\epf

\section{Formulas involving $\epsilon_0$} \lab{formulas}
 
For $S = \Spec \bC$, having proved the existence and uniqueness of $\eps_0$ in \cite{Del73}*{\S4}, Deligne proceeds to establish a formulary \cite{Del73}*{\S5} that details its properties and facilitates its computation. We gather some of these formulas here with a twofold aim: their special cases will be used in the proof of \Cref{main} to argue passage to more general bases $S$, and with little additional effort we will establish \eqref{ch-psi}, \eqref{unr-twist}, and \eqref{unr-twist-W} for all normal integral $S$.

\bpp[Change of additive character] \lab{f1}
As we have already observed in \S\ref{add-char}, $K^\times$ acts on the set of possible choices of $\psi$. The effect that this action bears on $\eps_0$ is explicated by
\be\lab{ch-psi}
\eps_0(V, a\psi, dx) = (\det V)(a) \cdot \abs{a}_K^{-\rk V} \cdot \eps_0(V, \psi, dx),
\ee
where $(\det V)(a)$ is the element of $\Gamma(S, \cO_S^\times)$ by which $a$ acts on the line bundle $\bigwedge^{\rk V} V$.
\epp

\bpp[Unramified twists] \lab{f2}
Twisting $V$ by an unramified $1$-dimensional character $\theta$ or, more generally, tensoring by an unramified $W$ of arbitrary dimension changes $\eps_0$ as follows:
\be\lab{unr-twist}
\eps_0(V\theta, \psi, dx) = \theta(\Frob_K)^{\Sw V + \rk V \cdot (n(\psi) + 1)} \eps_0(V, \psi, dx),
\ee
\be\lab{unr-twist-W}
\eps_0(V\tensor W, \psi, dx) = (\det W)(\Frob_K)^{\Sw V + \rk V \cdot (n(\psi) + 1)} \eps_0(V, \psi, dx)^{\rk W}.
\ee
\epp

\bpp[Explicit inverse in the $1$-dimensional case] \lab{1-dim-expl}
Unlike the other formulas, we will deduce \eqref{expl-inv} only for spectra of fields $F$ of characteristic different from $p$ from the known $F = \bC$ case\footnote{One difficulty encountered over more general bases is the incompatibility of $\rk \chi^I$ with reduction modulo $l$.}, and this will be one of the key inputs in proving that the $\eps_0$ are global units as claimed in \Cref{main}.

Given $\psi$ and $dx$, the \emph{dual Haar measure} $\wh{dx}$ of $dx$ with respect to $\psi$ is defined by insisting that 
\[
\int_{\cO_K} dx \cdot \int_{\cO_K} \wh{dx} = (\#\bF_K)^{-n(\psi)}.
\]
The resulting $\wh{dx}$ is well-defined thanks to the discussion of \S\ref{haar}.

Suppose that $S = \Spec \bC$, and let $C^\infty_c(K)$ be the space of locally constant compactly supported $\bC$-valued functions on $K$. Then the composition of the Fourier transform on $C^\infty_c(K)$ with respect to $\psi$ and $dx$ with the Fourier transform with respect to\footnote{Here $-\psi$ should be interpreted as $(-1) \cdot \psi$ with $-1 \in K$, see \S\ref{add-char}.} $-\psi$ and $\wh{dx}$ is the identity: compare, e.g., \cite{BH06}*{\S23.1, proof of Prop.}\footnote{Beware that what loc.~cit.~calls the level of $\psi$ is $-n(\psi)$ in the notation used here.}. Therefore, Tate's local functional equation \cite{Del73}*{5.8.1} for a continuous $\chi \in \Hom(K^\times, \bC^\times)$ gives
\[
\eps(\chi, \psi, dx)\eps(\chi\i \abs{\cdot}_K, -\psi, \wh{dx}) = 1,
\]
where $\eps(\chi, \psi, dx) = -\chi(\Frob_K)^{-\rk \chi^I} \eps_0(\chi, \psi, dx)$ as in Remark \ref{full-eps}. This proves the $F = \bC$ case of
\be\lab{expl-inv}
\eps_0(\chi, \psi, dx)\eps_0(\chi\i \abs{\cdot}_K, -\psi, \wh{dx}) = (\#\bF_K)^{-\rk \chi^I},
\ee
which is often helpful when reasoning about the inverse of $\eps_0(\chi, \psi, dx)$.
\epp

\bpp[Explicit inverse in the higher-dimensional case] \lab{high-dim}
Assume that $S = \Spec F$ for a field $F$ of characteristic $l$ with $l \neq p$. If $l = 0$, then set $I\pr \ce I$. If $l > 0$, then set $I\pr$ to be the preimage in $I$ of the compositum of the prime-to-$l$ Sylow subgroups of the quotient of $I$ by the wild inertia. In both cases one may interpret $I\pr$ as the minimal subgroup of $I$ for which $I/I\pr$ is pro-$l$. The closed subgroup $I\pr$ is normal in $W(K^s/K)$. Moreover, the finite quotients of $I\pr$ are of order prime to $l$, so the functor of taking $I\pr$-invariants is exact. Therefore, this functor induces a homomorphism
\[
(-)^{I\pr}\colon R_F(W(K^s/K)) \ra R_F(W(K^s/K)/I\pr).
\]
Let $V$ be a continuous representation of $W(K^s/K)$ over $F$. If $\dim V = 1$, then the action of $I$ on $V$ factors through a finite quotient of prime to $l$ order; thus, $V^{I\pr} = V^I$ for such $V$. Consequently,
\be \lab{expl-inv-2}
\eps_0(V, \psi, dx)\eps_0(V^* \abs{\cdot}_K, -\psi, \wh{dx}) = (\#\bF_K)^{-\rk V^{I\pr}},
\ee
where $V^*$ denotes the dual representation $\Hom_F(V, F)$, is an extension of \eqref{expl-inv} beyond the $1$-dimensional case. We will prove \eqref{expl-inv-2} for all $F$ at once in \Cref{fields-expl-inv}. The formula \eqref{expl-inv-2} will not play a role in the proof of \Cref{main}.
\epp

\section{The case when $S$ is restricted to spectra of fields of characteristic $0$} \lab{char-0}

As remarked in \ref{hist} and \ref{S=C}, for $S = \Spec \bC$ there exists a unique $\eps_0$ satisfying \ref{add}--\ref{dim-1}. Consequently,
\be\lab{aut-C}\tag{$\dagger$}
\text{\ref{BC} holds whenever }f\text{ arises from an element of }\Aut(\bC).
\ee

\bprop \lab{char-0-concl}
\Cref{main} holds if one restricts to $S$ of the form $\Spec F$ for a field $F$ of characteristic $0$. Moreover, the resulting $\eps_0$ satisfies \eqref{ch-psi}, \eqref{unr-twist}, \eqref{unr-twist-W}, and \eqref{expl-inv}.
\eprop

\bpf
Since $I_K$ acts through a finite quotient, $V \cong V\pr\tensor_{F\pr} F$ for a subfield $F\pr \subset F$ of finite transcendence degree over $\bQ$ and a representation $V\pr$ of $W(K^s/K)$ over $F\pr$. Enlarging $F\pr$ if needed, we assume further that $\psi$ and the Haar measure are $F\pr$-valued. Due to \ref{BC}, a choice of an embedding $\iota\colon F\pr \hra \bC$ forces us to set 
\[
\eps_0(V, \psi, dx) \ce \iota\i(\eps_0(V\pr \tensor_{F\pr, \iota} \bC, \iota\circ \psi, \iota \circ dx)).
\]
Once we check, as we do below, that the resulting $\eps_0(V, \psi, dx)$ is independent of choices, \ref{BC}--\ref{dim-1} as well as the claimed formulas will follow from the construction and the assumed $F = \bC$ case.

Firstly, $\eps_0(V\pr \tensor_{F\pr, \iota} \bC, \iota\circ \psi, \iota \circ dx) \in \iota(F\pr)^\times$ due to \eqref{aut-C}, because $(\bC^\times)^{\Aut(\bC/\iota(F\pr))} = \iota(F\pr)^\times$ \cite{BouA}*{V.107, Prop.~10}. Moreover, $\Aut(\bC)$ acts transitively on the set of embeddings of $F\pr$ into $\bC$ \cite{BouA}*{V.107, Cor.~2}, so \eqref{aut-C} also shows the independence of $\eps_0(V, \psi, dx)$ of the choice of $\iota$, and hence the independence of enlarging $F\pr$, as well. The independence of the choice of $V\pr$ follows, too, because any two choices are isomorphic over a larger $F\pr$.
\epf

\section{The case when $S$ is restricted to spectra of fields} \lab{field-case}

To settle this case in stages in \Cref{alg-cl,fields}, we discuss the necessary representation-theoretic preliminaries in \ref{elem}--\ref{R0-alg-cl}. To prepare for those, we recall the following well known lemma.

\blem \lab{dvr-res}
Fix a finite group $G$, and let $F$ be a sufficiently large field. There exists a complete discrete valuation ring $A$ with the residue field $F = A/\fm$ and the field of fractions $\eta$ that is sufficiently large and of characteristic $0$.
\elem

\bpf
For $F$ of characteristic $0$, one takes $A = F\llb t\rrb$. For $F$ of characteristic $l > 0$, one lets $m$ be the least common multiple of the orders of elements of $G$, applies \cite{Mat89}*{29.1} to $\bQ_l(\zeta_m)$, and replaces the resulting discrete valuation ring $A$ by its completion.
\epf

\bpp[Elementary groups] \lab{elem}
For a prime $p$, a finite group $H$ is \emph{$p$-elementary} if it is a product of a $p$-group and a cyclic group of order prime to $p$. A finite group is \emph{elementary} if it is $p$-elementary for some prime $p$. An elementary $H$ is the direct product of its Sylow subgroups; consequently, every $h \in H$ is a product of commuting elements of prime power order that are powers of $h$. To conclude that every $H\pr \le H$ is again the direct product of its Sylow subgroups, and hence also elementary, it remains to note that the latter are the intersections of $H\pr$ with the Sylow subgroups of $H$.
\epp

Many subsequent arguments will be based on the following version of Brauer's induction theorem.

\bprop[\cite{Del73}*{1.5}] \lab{brau}
For a finite group $G$ and a sufficiently large field $F$,
\benum
\item \lab{brau-a}
$R_F(G)$ is spanned by the elements of the form $\Ind_H^G [\chi]$ for elementary $H \le G$ and characters $\chi \in \Hom(H, F^\times)$, and

\item \lab{brau-b}
$R_F^0(G)$ is spanned by the elements of the form $\Ind_H^G([\chi] - [\mathbf{1}_H])$ for elementary $H \le G$ and characters $\chi \in \Hom(H, F^\times)$.
\eenum
\eprop

\bpf
Suppose initially that $\Char F = 0$. By Brauer's induction theorem \cite{Ser77}*{\S12.6, Thm.~27},
\[
[\b{1}_G] = \textstyle{\sum_i} \Ind_{H_i}^G a_i\quad\text{ in }R_F(G) \text{ for some elementary } H_i\le G \text{ and } a_i \in R_F(H_i).
\]
After multiplying both sides of this equality by a $v \in R_F(G)$ (resp.,~a $v \in R_F^0(G)$ for \ref{brau-b}), the projection formula reduces to the case when $G$ is elementary itself. On the other hand, by \cite{Ser77}*{\S12.3, proof of Thm.~24}, $R_F(G)$ is spanned by the elements of the form $\Ind_H^G[\chi]$ for subgroups $H \le G$ and characters $\chi \in \Hom(H, F^\times)$, so \ref{brau-a} follows because, as noted in \S\ref{elem}, $H$ is elementary if so is $G$. As for \ref{brau-b}, this expresses every $v \in R_F^0(G)$ as 
\[
v = \textstyle{\sum_i} n_i\Ind_{H_i}^G [\chi_i] = \textstyle{\sum_i} n_i \Ind_{H_i}^G ([\chi_i] - [\b{1}_{H_i}]) + \textstyle{\sum_i} n_i \Ind_{H_i}^G [\b{1}_{H_i}],\quad\quad n_i \in \bZ
\]
for elementary $H_i \le G$ and characters $\chi_i \in \Hom(H_i, F^\times)$, so the conclusion follows by induction on $\#G$ because $\textstyle{\sum_i} n_i \Ind_{H_i}^G [\b{1}_{H_i}] \in R_F^0(G)$ is the inflation of an element of $R_F^0(G/Z)$ where $Z \normal\, G$ is the center, which is nontrivial if so is the elementary $G$.

Suppose now that $F$ is of characteristic $l > 0$ and choose $A$ as in \Cref{dvr-res}. As observed in \S\ref{dec}, the decomposition homomorphism $d_G\colon R_\eta(G) \ra R_F(G)$ is surjective, preserves ranks, and commutes with induction, so we deduce \ref{brau-a} and \ref{brau-b} for $F$ from the case of $\eta$ established above.
\epf

\bprop[\cite{Del73}*{1.8}] \lab{dec-ker}
For a finite group $G$ and a discrete valuation ring $A$ with the residue field $F = A/\fm$ of characteristic $l$ and the field of fractions $\eta$ that is sufficiently large and of characteristic $0$, the elements of the form $\Ind_H^G([\chi] - [\chi\pr])$ for elementary $H \le G$ and characters $\chi, \chi\pr \in \Hom(H, A^\times)$ with $\chi(h) \equiv \chi\pr(h) \bmod \fm$ for all $h \in H$ generate $\Ker(d_G\colon R_\eta(G) \ra R_F(G))$.
\eprop

\bpf 
Passing to $\wh{\eta}$ as in \S\ref{dec}, one may for comfort reduce to the complete case. Also, $d_G$ is an isomorphism if $l = 0$ (see \S\ref{dec}), so we assume for the remainder of the proof that $l > 0$.

Using \Cref{brau} \ref{brau-a} to write 
\[
[\b{1}_G] = \textstyle{\sum_i} n_i \Ind_{H_i}^G [\chi_i] \quad \text{ in $R_\eta(G)$ for some $n_i \in \bZ$, elementary $H_i \le G$, and $\chi_i \in \Hom(H_i, \eta^\times)$,}
\]
for $v \in \Ker d_G$ we have
\[
v = \textstyle{\sum_i} n_i \Ind_{H_i}^G( [\chi_i] \cdot \Res_{H_i}^G v).
\]
Since $\Res_{H_i}^G v \in \Ker d_{H_i}$, we are reduced to the case of an elementary $G$.

An elementary $G$ is a product $G = N \times P$ with $l \nmid \#N$ and $\#P = l^n$ for some $n \ge 0$. By \cite{CR81}*{10.33 and 17.1}, the irreducible representations of $G$ over $\eta$ or $F$ are precisely the tensor products of an irreducible representation of $N$ with a one of $P$. Consequently, $d_G$ takes the form
\[
R_\eta(G) \cong R_\eta(N) \tensor_\bZ R_\eta(P) \xra{d_N \tensor d_P} R_F(N) \tensor_\bZ R_F(P) \cong R_F(G).
\]
Since $d_N$ is an isomorphism \cite{Ser77}*{\S15.5} and $R_\eta(N)$ is $\bZ$-free, 
\[
\Ker d_G \cong R_\eta(N) \tensor \Ker d_P,
\]
so \Cref{brau} \ref{brau-a} applied to $N$ reduces further to the case $G = P$.

However, if $G$ is of $l$-power order, then every character $\chi \in \Hom(H, \eta^\times)$ of a subgroup $H \le G$ takes values in $1 + \fm A \subset A^\times$, and the claim results from \Cref{brau} \ref{brau-b}.
\epf

\bpp[Weil representation types (\cite{Del73}*{4.10})] \lab{types} 
Fix a separably closed field $F$, and let $l$ be its characteristic. A continuous representation $V$ of $W(K^s/K)$ over $F$ is said to have a type if $\Frob_K^{m}$ acts on $V$ as a scalar $a \in F^\times$ for some $m \ge 1$. When this is the case, $\Frob_K^{Nm}$ acts on $V$ as the scalar $a^{N}$ for $N \ge 1$, so the resulting element of the direct limit (which is indexed by the positive integers ordered by the divisibility relation)
\[
\varinjlim_{m \mid n} F^\times\quad \text{ with transition maps }\quad a \mapsto a^{\f{n}{m}}\quad \text{ between the copies of }F^\times\text{ in positions }m\text{ and }n
\]
 is independent of the choice of $m$; it is the \emph{type} of $V$. To argue that the type is also independent of the choice of $\Frob_K$, let $W$ be the quotient through which $W(K^s/K)$ acts on $V$, let $J \normal W$ be the finite image of $I_K$, and write $\ov{\Frob}_K \in W$ for the image of $\Frob_K$. Since $J$ is centralized by $\ov{\Frob}_K^m$ for all sufficiently divisible $m$ and changing $\Frob_K$ changes $\ov{\Frob}_K^m$ by an element of $J$, we see that $\ov{\Frob}_K^{m \cdot \#J}$ is independent of $\Frob_K$, and hence so is the type.
\epp

\begin{lem} \lab{irr-type}
Every irreducible $V$ has a type.
\end{lem}

\bpf
Indeed, $\End_{W} V$ is a finite dimensional division algebra over the separably closed $F$, so $\End_{W} V \cong F\pr$ for a finite field extension $F\pr/F$. Thus, every $\ov{\Frob}_K^N$ that is central in $W$ acts as a scalar in $F^{\prime\times}$, and hence $\ov{\Frob}_K^{Nl^a}$ acts as a scalar in $F^\times$ for large $a$. 
\epf

For a fixed type $\tau$, let $R_{F, \tau}(W(K^s/K))$ be the subgroup of $R_F(W(K^s/K))$ spanned by the classes of representations of type $\tau$. Due to \Cref{irr-type} and the discussion of \S\ref{Gr-fields},
\be\lab{dec-types}
R_F(W(K^s/K)) = \textstyle\bigoplus_{\tau} R_{F, \tau}(W(K^s/K)).
\ee
The representations of type $1$ are precisely the Galois representations, see \S\ref{weil}. 

If $F$ is in addition algebraically closed, then for every type $\tau$ there is an unramified character $\chi_\tau\colon W(K^s/K) \ra F^\times$ of this type. Twisting by $\chi_\tau$ induces an isomorphism 
\be\lab{twist-types}
\xymatrix{
R_{F, 1}(W(K^s/K)) \ar[rr]^{-\tensor \chi_\tau}_{\sim} &&R_{F, \tau}(W(K^s/K)).
}
\ee

\bprop \lab{R0-alg-cl}
For an algebraically closed $F$, the elements of the form $\Ind_{W(K^s/L)}^{W(K^s/K)}([\chi] - [\chi\pr])$ for finite subextensions $K^s/L/K$ and continuous $\chi, \chi\pr \in \Hom(W(K^s/L), F^\times)$ span $R_F^0(W(K^s/K))$.
\eprop

\bpf
For a $v \in R_F^0(W(K^s/K))$, let $v = \sum_\tau v_\tau$ be the decomposition provided by \eqref{dec-types}. For each appearing $\tau$, fix a $\chi_\tau$ as in \eqref{twist-types}. Then  $\sum_\tau \rk v_\tau \cdot [\chi_\tau]$ is in the desired span ($L = K$ suffices), and it remains to note that so is each $v_\tau - \rk v_\tau \cdot [\chi_\tau]$ thanks to \Cref{brau} \ref{brau-b} and \eqref{twist-types}.
\epf

Having settled the representation-theoretic preliminaries, we get back to studying $\eps_0$.

\blem \lab{ker-eps}
Let $A$ be a discrete valuation ring with the residue field $F = A/\fm$ of characteristic~$l$ with $l \neq p$ and the fraction field $\eta$ of characteristic $0$. Adopt the setup of \Cref{main} for $S = \Spec A$. 
\benum
\item \lab{ker-eps-a}
For the $\eps_0$ of \Cref{char-0-concl}, one has $\eps_0(V_{\eta}, \psi_\eta, (dx)_\eta) \in A^\times$.

\item \lab{ker-eps-b}
For a finite Galois subextension $K^s/L/K$ with the Galois group $G = \Gal(L/K)$, the restriction of the $\eps_0$ of \Cref{char-0-concl} to the kernel of $R_\eta(G) \ra R_F(G)$ takes values in $1 + \fm A \subset \eta^\times$. 
\eenum
\elem

\bpf \hfill
\benum
\item
Due to \ref{BC} and the existence of a discrete valuation prolonging the given one on $\eta$ to any finite extension $\eta\pr/\eta$ \cite{Mat89}*{Cor.~to 11.7}, we may replace $\eta$ by any $\eta\pr$ if needed. Moreover, for a $W(K^s/K)$-stable $A$-lattice in $V_\eta$, its intersection with a subrepresentation (resp.~image in a quotient) is a stable $A$-lattice, so any subquotient of $V_\eta$ is realizable over $A$. Therefore, we may assume that $V_\eta$ is absolutely irreducible, there is an unramified character $\theta \colon W(K^s/K) \ra \eta^\times$ with the same type as $V_{\ov{\eta}}$, and the Galois representation $V_\eta \tensor \theta\i$ factors through a continuous finite quotient $G$ of $W(K^s/K)$ for which $\eta$ is sufficiently large. Since scaling by $\theta(\Frob_K^N)$ is an automorphism of the $A$-lattice $V \subset V_\eta$ if $N$ is sufficiently divisible, $\theta$ takes values in $A^\times$. Thus, \eqref{unr-twist}, \Cref{brau} \ref{brau-b} applied to $V_\eta \tensor \theta\i - \rk V \cdot [\b{1}_G]$, and \ref{ind} reduce to the $1$-dimensional case, in which the $A$-valued $\eps_0$ given by \ref{dim-1} has the $A$-valued inverse provided explicitly by \eqref{expl-inv} and \ref{dim-1}. 

\item
As in the proof of \ref{ker-eps-a}, we assume that $\eta$ is sufficiently large.
\Cref{dec-ker} then reduces to proving that $\eps_0(\Ind_H^G([\chi] - [\chi\pr]), \psi, dx) \in 1 + \fm A$ for every subgroup $H \le G$ and characters $\chi, \chi\pr \in \Hom(H, A^\times)$ with $\chi(h) \equiv \chi\pr(h) \bmod \fm$ for all $h \in H$. The additive character $\psi$ necessarily takes values in $A$. The property \ref{ind} reduces further to proving that
\be \lab{very-temp}
\quad \eps_0(\chi, \psi \circ \Tr_{L^H/K}, dx) \equiv \eps_0(\chi\pr, \psi \circ \Tr_{L^H/K}, dx) \bmod \fm \quad \quad \text{with both sides in $A$}
\ee
for a $dx$ that also takes values in $A$. To obtain \eqref{very-temp}, apply \ref{dim-1} and note that $\Sw(\chi) = \Sw(\chi\pr)$, as can be checked over the closed point thanks to \Cref{Swan-sch}. \qedhere
\eenum
\epf

\bprop \lab{alg-cl}
In \Cref{main}, restrict to $S$ of the form $\Spec F$ for an algebraically closed~field~$F$.
\benum

\item \lab{alg-cl-Gal}
An $\eps_0$ satisfying \ref{add}--\ref{dim-1}, \eqref{ch-psi}, \eqref{unr-twist}, and \eqref{unr-twist-W} exists if it does when $V$, $\theta$, and $W$ are restricted to Galois representations. 

\item  \lab{unique}
If one restricts $V$ to Galois representations, then an $\eps_0$ satisfying \ref{add}--\ref{dim-1} is unique if it exists, in which case it also satisfies \ref{BC}. The same conclusion holds if one restricts further to representations that factor through a fixed finite quotient $\Gal(L/K)$.

\item \lab{fields-alg}
There exists a unique $\eps_0$ satisfying \ref{BC}--\ref{dim-1}. It also satisfies \eqref{ch-psi}, \eqref{unr-twist}, and \eqref{unr-twist-W}.
\eenum
\eprop

\bpf
For every $A$ provided by \Cref{dvr-res}, Hensel's lemma (or, if one prefers, \cite{EGAIV4}*{18.5.15}) lifts every $\psi$ valued in $F^{\times}$ to a unique additive character valued in $A^\times$; lifting $dx$ amounts to lifting $\int_{\cO_K} dx \in F^\times$ to $A^\times$, see \S\ref{haar}. We continue to denote these lifts by $\psi$ and $dx$ (although the latter is not unique) and recall from \S\ref{add-char} that $n(\psi)$ is invariant under reduction modulo $\fm$. Lifting is also possible for unramified characters of Weil groups and $1$-dimensional characters of finite groups.

\benum
\item
For each type $\tau$, take an unramified character $\chi_\tau\colon W(K^s/K) \ra F^\times$ of this type. To define $\eps_0$ when its restriction to Galois representations $v \in R_{F, 1}(W(K^s/K))$ is given, set 
\be\lab{tw-def}
\eps_0(v \tensor \chi_\tau, \psi, dx) \ce \chi_\tau(\Frob_K)^{\Sw v + \rk v \cdot (n(\psi) + 1)} \cdot \eps_0(v, \psi, dx)
\ee
and extend $\eps_0$ to $R_F(W(K^s/K))$ using \eqref{dec-types} and \eqref{twist-types}. Due to \eqref{unr-twist} for Galois representations, the definition \eqref{tw-def} does not depend on the choice of $\chi_\tau$. 

The desired \ref{add}, \ref{haar-scale}, and \ref{dim-1} are immediate. So is \eqref{ch-psi}, because $n(a\psi) = n(\psi) + v_K(a)$ where $v_K(a)$ is the valuation of $a$. %Once \ref{ind} is proved, \eqref{unr-twist} will follow from \ref{unique}. 
In checking \eqref{unr-twist-W}, which includes \eqref{unr-twist}, one restricts to irreducible $V$ and $W$, which share their types with unramified characters $\chi_V$ and $\chi_W$. Writing $W \cong W\pr \tensor \chi_W$ and $V \cong V\pr \tensor \chi_V$, \eqref{unr-twist-W} results from its version for $V\pr \tensor W\pr$, \eqref{tw-def} with $\chi_\tau = \chi_V\chi_W$, and the relation $\det W = \chi_W^{\rk W} \det W\pr$.

We concentrate on the remaining \ref{ind}; it can be proved by a small computation, which, however, can also be relegated to characteristic $0$, as we now explain. \Cref{R0-alg-cl} reduces to considering $v \in R_F^0(W(K^s/L))$ of the form $[\chi] - [\chi\pr]$ for continuous $\chi, \chi\pr \in \Hom(W(K^s/L), F^\times)$. For such a $v$, choose unramified $\theta, \theta\pr \in \Hom(W(K^s/K), F^\times)$ for which $\theta|_L$ and $\theta\pr|_L$ have the same types as $\chi$ and $\chi\pr$. Let $G$ be a continuous finite quotient of $W(K^s/L)$ through which $\chi\i\cdot \theta|_L$ and $\chi^{\prime -1}\cdot \theta\pr|_L$ factor, and take $A$ as in \Cref{dvr-res}. Lifting $\theta$, $\theta\pr$, $\chi\i\cdot \theta|_L$, $\chi^{\prime -1}\cdot \theta\pr|_L$, $\psi$, and the involved Haar measures to $A$ yields the desired conclusion thanks to \Cref{char-0-concl}, as long as we exhibit $\eps_0(\Ind_{W(K^s/L)}^{W(K^s/K)} \chi, \psi, dx_K)$, $\eps_0(\chi, \psi \circ \Tr_{L/K}, dx_L)$, and their $\chi\pr$ analogues  as reductions of the corresponding elements of $A^\times$ (cf.~\Cref{ker-eps} \ref{ker-eps-a}). 

To argue this, thanks to \Cref{Swan-sch} and the unramified twist aspect of \Cref{char-0-concl}, we only need to prove that the formation of $\eps_0$ is compatible with reduction modulo $\fm$ for Galois representations. Having the liberty of replacing $A$ by its normalization in a finite extension of $\eta$, we deduce this compatibility from \ref{ind}, \ref{dim-1}, and \Cref{brau} \ref{brau-b}.

\item
To show that the value of $\eps_0$ on $V$ is uniquely determined, take a finite quotient $\Gal(L/K)$ through which the $W(K^s/K)$-action factors and use \Cref{brau} \ref{brau-b} to express
\[
[V] - \rk V \cdot [\b{1}_K] = \sum \Ind_{W(K^s/\wt{L})}^{W(K^s/K)}\p{[\chi] - [\b{1}_{\wt{L}}]}\quad \text{in $R_F(W(K^s/K))$}
\]
for suitable subextensions $L/\wt{L}/K$ and characters $\chi \in \Hom(\Gal(L/\wt{L}), F^\times)$. The uniqueness, as well as \ref{BC}, follows immediately from \ref{add}--\ref{dim-1}.

\item
Thanks to \ref{unique} and \ref{alg-cl-Gal}, we seek an $\eps_0$ satisfying \ref{add}--\ref{dim-1}, \eqref{ch-psi}, \eqref{unr-twist}, and \eqref{unr-twist-W} for Galois representations. Furthermore, we restrict to representations that factor through a fixed finite quotient $G = \Gal(L/K)$, since the resulting $\eps_0$ for varying $G$ will be compatible by \ref{unique}. 

For an $A$ given by \Cref{dvr-res} and the $\eps_0$ constructed in \Cref{char-0-concl}, \Cref{ker-eps} \ref{ker-eps-a} gives
\[ 
\eps_0(-, \psi, dx)\colon R_\eta(G) \ra A^\times \subset \eta^\times.
\] 
Thus, due to \Cref{ker-eps} \ref{ker-eps-b} and the surjectivity of $d_G$ (see \S\ref{dec}), $\eps_0(-, \psi, dx) \bmod \fm$ induces
\[
\eps_0(-, \psi, dx)\colon R_F(G) \ra F^\times,
\]
which satisfies \ref{add}--\ref{dim-1}, \eqref{ch-psi}, \eqref{unr-twist}, and \eqref{unr-twist-W} by construction and \Cref{Swan-sch}. 
\qedhere
\eenum
\epf

\bprop \lab{fields}
\Cref{main} holds if one restricts to $S$ of the form $\Spec F$ for a field $F$. Moreover, the resulting $\eps_0$ satisfies \eqref{ch-psi}, \eqref{unr-twist}, and \eqref{unr-twist-W}.
\eprop

We will show in \Cref{fields-expl-inv} that $\eps_0$ also satisfies \eqref{expl-inv} and \eqref{expl-inv-2}.

\bpf
Since \ref{BC} forces $\eps_0(V, \psi, dx) \ce \eps_0(V_{\ov{F}}, \psi, dx) \in \ov{F}^\times$, thanks to \Cref{alg-cl} \ref{fields-alg}, we only need to check that the resulting $\eps_0$ takes values in $F$. Also, \ref{BC} applied to $f$ arising from the elements of $\Aut(\ov{F}/F)$ reduces to the case of a separably closed $F$. We therefore assume that $F$ is separably closed and imperfect of characteristic $l > 0$. We may and do further assume that $V$ is irreducible and, scaling $dx$ if needed, that $\psi$ and $dx$ take values in $\ov{\bF}_l \subset F$. If $V$ is a Galois representation, then the claim results from \Cref{brau} \ref{brau-b}, \ref{ind}, and \ref{dim-1}.

Since $I_K$ is a normal subgroup, every element of $W(K^s/K)$ maps every $I_K$-stable subspace of $V$ to another such subspace. Consequently, the sum of the irreducible $I_K$-subrepresentations of $V$ is $W(K^s/K)$-stable, so the irreducibility of $V$ entails the semisimplicity of $V|_{I_K}$. Moreover, $W(K^s/K)$ acts transitively on the $I_K$-isotypic components $V\pr \subset V$, so $V \cong \Ind_{W(K^s/L)}^{W(K^s/K)} V\pr$ for such a $V\pr$ and its stabilizer $W(K^s/L) \le W(K^s/K)$ (compare \cite{Ser77}*{\S8.1, proof of Prop.~24} if needed). Therefore, since $\Ind_{W(K^s/L)}^{W(K^s/K)} [\b{1}_L]$ is a Galois representation, induction on $\dim V$ by means of \ref{ind} allows us to assume that $V|_{I_K}$ is isotypic. % If $V\pr$ is not isotypic for the inertia of $W(K^s/L)$, then at least its dimension is smaller and one can repeat the process.

Using \S\ref{Gr-fields}, we let $(\rho, X)$ be the irreducible representation of $I_K$ over $\ov{\bF}_{l}$ for which $V|_{I_K}$ is a multiple of $X_F$. Since the isomorphism class of $X_F$ is preserved under $\Frob_K$-conjugation, so is that of $X$, see \S\ref{Gr-fields}. A choice of an $\ov{\bF}_{l}$-isomorphism $\iota\colon (\rho, X) \isomto (\rho\pr, X)$ where $\rho\pr(i) = \rho(\Frob_K i\Frob_K\i)$ for $i \in I_K$ extends $X$ to a representation of $W(K^s/K)$ over $\ov{\bF}_{l}$ by letting $\Frob_K$ act as $\iota$. This extension, still denoted by $X$, is a Galois representation because $\iota \in \GL(X)$ has finite order. 

Since $\dim_F \End_{I_K}(X_F, X_F) = \dim_{\ov{\bF}_l} \End_{I_K}(X, X) = 1$ (for the last equality, see the proof of \Cref{irr-type}), $\End_{I_K}(X_F, X_F) \cong F$, so the canonical $X_F \tensor_F \Hom_{I_K}(X_F, V) \ra V$ is an isomorphism. Evidently, it is also $W(K^s/K)$-equivariant, so $V$ decomposes over $F$ as a tensor product of the Galois representation $X_F$ and the unramified $\Hom_{I_K}(X_F, V)$. It remains to apply \eqref{unr-twist-W}.
\epf

We record the following strengthening of \Cref{ker-eps} \ref{ker-eps-a} that will be used in \S\ref{genl}.

\bprop \lab{red-dvr}
Let $A$ be a discrete valuation ring with the residue field $F$ of characteristic $l$ with $l \neq p$ and the fraction field $\eta$. Adopt the setup of \Cref{main} for $S = \Spec A$. For the $\eps_0$ of \Cref{fields}, one has $\eps_0(V_\eta, \psi_\eta, (dx)_\eta) \in A^\times$ with the image $\eps_0(V_F, \psi_F, (dx)_F)$ in $F^\times$. 
\eprop

\bpf
The proof is similar to that of \Cref{ker-eps} \ref{ker-eps-a}. Namely, if needed we again replace $\eta$ by a finite extension and $V$ by its subquotient to assume that $V \cong V\pr \tensor \chi_V$ for an $A^\times$-valued unramified character $\chi_V$ and a Galois representation $V\pr$ that factors through a continuous finite quotient $G$ of $W(K^s/K)$ for which $\eta$ is sufficiently large. The formula \eqref{unr-twist} (together with \S\ref{add-char} and \Cref{Swan-sch}) then reduces to the case $V = V\pr$, and we can further assume that $\dim V = 1$ thanks to \Cref{brau} \ref{brau-b}, \ref{ind}, and the compatibility of the decomposition homomorphism with induction (see \S\ref{dec}). However, in the $1$-dimensional Galois case, $\eps_0(V_\eta, \psi_\eta, (dx)_\eta) \in A$ thanks to \eqref{dim-1-def}, which also shows that this $\eps_0$ reduces to $\eps_0(V_F, \psi_F, (dx)_F)$ in $F$. Since the latter is nonzero by \Cref{fields}, the claim follows.
\epf

\bcor \lab{fields-expl-inv}
The $\eps_0$ of \Cref{fields} also satisfies \eqref{expl-inv} and \eqref{expl-inv-2}.
\ecor

\bpf
We begin with \eqref{expl-inv}. Let $\chi_0 \in \Hom(W(K^s/K), F^\times)$ be the unramified character for which $\chi_0(\Frob_K) = \chi(\Frob_K)$. Replacing $\chi$ by $\chi\chi_0\i$, invoking \eqref{unr-twist}, and passing to $\ov{F}$, we may and do assume that $\chi$ factors through a continuous finite quotient $G$ for which $F$ is sufficiently large. We then use the remarks in the beginning of the proof of \Cref{alg-cl} to lift $\psi$, $dx$, and $\chi$ to corresponding objects over an $A$ provided by \Cref{dvr-res}. Since the lift of $\chi$ is ramified if and only if so is $\chi$, \Cref{red-dvr} reduces to the $\Char F = 0$ case, which was established in \Cref{char-0-concl}.

We now turn to \eqref{expl-inv-2}, for which we set $l \ce \Char F$ and also adopt other notation of \S\ref{high-dim}. Both sides of \eqref{expl-inv-2} define homomorphisms $R_F(W(K^s/K))\ra F^\times$: the left one due to \ref{add} and the exactness of $V \mapsto V^*$, the right one due to the exactness of $(-)^{I\pr}$ noted in \S\ref{high-dim}. Moreover, by \eqref{expl-inv}, the two homomorphisms agree on the classes of $1$-dimensional representations. Passing to $\ov{F}$ and invoking \Cref{R0-alg-cl}, we therefore reduce to showing that they also agree on the elements of the form $\Ind_{W(K^s/L)}^{W(K^s/K)}([\chi] - [\chi\pr])$ for finite subextensions $K^s/L/K$ and continuous $\chi, \chi\pr \in \Hom(W(K^s/L), F^\times)$. Taking duals commutes with induction,\footnote{This commutation follows from the self-duality of $F[W(K^s/K)/W(K^s/L)]$, which in turn follows from the $W(K^s/K)$-invariance of the nondegenerate bilinear pairing $F[W(K^s/K)/W(K^s/L)] \times F[W(K^s/K)/W(K^s/L)] \ra F$ for which the standard basis (indexed by the cosets) is its own dual basis.} so the combination of \ref{haar-scale}, \ref{ind}, and \eqref{expl-inv} shows that the left hand side of \eqref{expl-inv-2} evaluates on $\Ind_{W(K^s/L)}^{W(K^s/K)}([\chi] - [\chi\pr])$ to
\[
\eps_0\p{[\chi] - [\chi\pr], \psi \circ \Tr_{L/K}}\eps_0\p{[\chi\i\abs{\cdot}_L] - [\chi^{\prime -1} \abs{\cdot}_L], -\psi \circ \Tr_{L/K}} = (\#\bF_L)^{-\rk \chi^{I\pr_L} + \rk \chi^{\prime I\pr_L}},
\]
where $I_L\pr \ce I\pr \cap W(K^s/L)$, which equals the $I\pr$ of $W(K^s/L)$. This agrees with what the right hand side evaluates to thanks to the following \Cref{cC} \ref{claim-b}.

\begin{claim} \lab{cC} For a representation $V$ of $W(K^s/L)$ over $F$, one has
\benum
\item \lab{claim-a}
$\p{\Ind_{W(K^s/L)}^{W(K^s/K)} V}^{I\pr} \cong \Ind_{W(K^s/L)/I_L\pr}^{W(K^s/K)/I\pr} V^{I_L\pr}$;

\item \lab{claim-b}
$(\#\bF_K)^{\rk \p{\Ind_{W(K^s/L)}^{W(K^s/K)} V}^{I\pr}} = (\#\bF_L)^{\rk V^{I_L\pr}}$ in $F$.
\eenum
\end{claim}

\bpf \hfill
\benum
\item
The argument below imitates an argument in \cite{Del73}*{proof of 3.8}. 

After consulting \cite{BH06}*{\S2.5, proof of Lemma} to reconcile with the definition of induction given in \ref{ind}, we may identify 
\[
\Ind_{W(K^s/L)}^{W(K^s/K)} V = \{ f\colon W(K^s/K) \ra V \ | \  f(yx) = y f(x) \ \text{ for }\ y \in W(K^s/L) \},
\]
 with $w \in W(K^s/K)$ acting by $wf\colon x \mapsto f(xw)$. Likewise,
\[
\quad \quad \Ind_{W(K^s/L)/I_L\pr}^{W(K^s/K)/I\pr} V^{I_L\pr} = \{ \ov{f}\colon W(K^s/K)/I\pr \ra V^{I_L\pr} \ | \  \ov{f}(yx) = y \ov{f}(x) \ \text{ for }\ y \in W(K^s/L)/I_L\pr \},
\]
 with $w \in W(K^s/K)/I\pr$ acting by $w\ov{f}\colon x \mapsto \ov{f}(xw)$. 
 
Due to these descriptions, $\Ind_{W(K^s/L)/I_L\pr}^{W(K^s/K)/I\pr} V^{I_L\pr} \subset \p{\Ind_{W(K^s/L)}^{W(K^s/K)} V}^{I\pr}$. For the converse inclusion, an $f \in \p{\Ind_{W(K^s/L)}^{W(K^s/K)} V}^{I\pr}$ is an inflation of an $\ov{f}\colon W(K^s/K)/I\pr \ra V$ that must take values in $V^{I_L\pr}$, because $\ov{f}(x) = \ov{f}(yx) = y\ov{f}(x)$ for $y \in I_L\pr$.

\item
We have $[W(K^s/K)/I\pr : W(K^s/L)/I_L\pr] = [W(K^s/K) : W(K^s/L)I\pr]$, which, up to a power of $l$ if $l > 0$, equals $[\bF_L : \bF_K]$. The combination of this and \ref{claim-a} gives the claim because the power of $l$ does not matter: $(\#\bF_K)^l = \#\bF_K$ in $F$ if $l > 0$. \qedhere
\eenum
\epf 

The agreement of the two homomorphisms $R_F(W(K^s/K))\ra F^\times$ establishes \eqref{expl-inv-2}.
\epf

\section{The general case} \lab{genl}

\bpp[Normal integral schemes]
Recall from \cite{EGAI}*{\S0, 4.1.4} or \cite{EGAIV2}*{5.13.5} that a scheme $S$ is \emph{normal} if its local rings are integrally closed domains. For an integral $S$, normality is equivalent to $\cO_S(U)$ being an integrally closed domain for every open $U \subset S$ \cite{EGAII}*{8.8.6.1}.
\epp

\bpp[Universally Japanese rings] \lab{japan}
For an episodic appearance below, recall from \cite{EGAIV2}*{7.7.1} that a ring $R$ is \emph{universally Japanese} if for every domain $R\pr$ that is a finite type $R$-algebra, the integral closure of $R\pr$ in $\Frac R\pr$ is a finite $R\pr$-module. For our purposes it suffices to know that every Dedekind domain whose fraction field has characteristic $0$ is universally Japanese \cite{EGAIV2}*{7.7.4}.
\epp

\begin{proof}[Proof of \Cref{main} and the formulas \eqref{ch-psi}, \eqref{unr-twist}, and \eqref{unr-twist-W}]
\Cref{fields} and \ref{BC} force us to define
\[
\eps_0(V, \psi, dx) \ce \eps_0(V_\eta, \psi_\eta, (dx)_\eta),
\]
where $\eta$ is the generic point of $S$. Once we check, as we do below, that the resulting $\eps_0$ is $\Gamma(S, \cO_S^\times)$-valued and its image in $k(s)^\times$ is $\eps_0(V_{s}, \psi_{s}, (dx)_{s})$ for every $s \in S$, \ref{BC}--\ref{dim-1}, as well as \eqref{ch-psi}, \eqref{unr-twist}, and \eqref{unr-twist-W}, will follow from \Cref{fields}. 

The promised checking can be done locally on $S$, so we assume that $S = \Spec A$ is affine and $V$ is free. Once we express $A$ as a filtered direct limit of Noetherian normal $\bZ[\f{1}{p}]$-subalgebras containing the values of $\psi$ and $dx$, since $W(K^s/K)$ acts on $V$ through a finitely generated quotient, limit arguments and the base change aspect of \Cref{fields} will descend $V$ to such a subalgebra and permit to assume further that $A$ is Noetherian. For such an expression, to capture the values of $\psi$, endow $A$ with a structure of an $R$-algebra, where $R = \bZ[\f{1}{p}][\zeta_{p\I}]$ if $\Char K = 0$ and $R = \bZ[\f{1}{p}][\zeta_p]$ if $\Char K = p$; then $\psi$ is the pullback of an $R^\times$-valued additive character. The normal $1$-dimensional $R\subset \bZ[\f{1}{p}][\zeta_{p\I}]$ is a Dedekind domain because it is Noetherian due to Cohen's theorem: every finite prime $l$ different from $p$ is finitely decomposed and unramified in $\bQ(\zeta_{p\I})$, so every prime $\fp \subset R$ is generated by the finitely generated $\fp \cap \bZ[\f{1}{p}][\zeta_{p^n}]$ for some $n \ge 0$. Consequently, as recalled in \S\ref{japan}, $R$ is universally Japanese, so finite type $R$-subalgebras of $A$ have Noetherian normalizations in their fraction fields and contain the values of $dx$ whenever they contain $\int_{\cO_K} dx$ (see \S\ref{haar}), which completes the reduction to the Noetherian case.

In the Noetherian case, $A$ is an intersection of discrete valuation rings \cite{Mat89}*{12.3 and 12.4 (i)}:
\[
A = \bigcap_{\Ht(\fp) = 1} A_\fp\quad \text{inside $\Frac A$,}
\]
for which \Cref{red-dvr} supplies the required checking, so the desired $\eps_0(V_\eta, \psi_\eta, (dx)_\eta) \in A^\times$ follows. Lastly, for a nongeneric $s\pr \in S$ the image of $\eps_0(V_\eta, \psi_\eta, (dx)_\eta)$ in $k(s\pr)^\times$  is $\eps_0(V_{s\pr}, \psi_{s\pr}, (dx)_{s\pr})$, as one sees by taking $s = \eta$ in \Cref{lem-EGA} and combining its conclusion with the base change aspect of \Cref{fields}.
\epf

\end{document}